\documentclass[12pt,twoside]{amsart}
\usepackage{amsmath, amsthm, amscd, amsfonts, amssymb, graphicx}
\usepackage[bookmarksnumbered, plainpages]{hyperref}

\newtheorem{thm}{Theorem}[section]

\newtheorem{lem}[thm]{Lemma}

\numberwithin{equation}{section}

\begin{document}

\title{\bf A Liouville theorem on complete non-K\"ahler manifolds}

\author{Yuang Li, Chuanjing Zhang and Xi Zhang}

\thanks{{\scriptsize
\hskip -0.4 true cm \textit{2010 Mathematics Subject Classification:}
53C55; 53C07; 58E20.
}
\newline The authors were supported in part by NSF in
China,  No. 11625106, 11571332 and 11721101.
}

\maketitle

\begin{abstract}
In this paper, we prove a Liouville theorem for holomorphic functions on a class of complete Gauduchon manifolds. This  generalizes a result
of Yau for complete K\"ahler manifolds to the complete non-K\"ahler case.
\end{abstract}

\vskip 0.2 true cm


\pagestyle{myheadings}
\markboth{\rightline {\scriptsize Y. Li et al.}}
         {\leftline{\scriptsize A Liouville theorem on complete non-K\"ahler manifolds}}

\bigskip
\bigskip


\section{ Introduction}

Let $(M, g)$ be a complete Riemannian manifold and $\Delta $ be the Beltrami-Laplacian of the Riemannian metric. In \cite{Y76}, Yau studied the equation
 \begin{equation}\label{log}
 \Delta \log u =f,
 \end{equation}
 and proved that if $0 <\int_{M} f dv_{g} \leq \infty $ or if $f \equiv 0$, then the  equation (\ref{log}) has no non-constant $L^{p}$-solutions for $0 < p <\infty$. As an application, Yau obtained the following $L^{p}$ Liouville theorem.

 \begin{thm}[\cite{Y76}, Theorem 4]\label{thmYau}
Let $(M, \omega )$ be a complete K\"ahler manifold. Then there is no non-constant $L^{p}$ holomorphic functions for $p>0$.
\end{thm}

It should be pointed out that, under the assumption that $(M, \omega )$ has sectional curvature of the same sign, Greene-Wu (\cite{GW71}) got some lower bound for the $L^{p}$ integral of  holomorphic functions. In this paper, we want to study the $L^{p}$
Liouville theorem of holomorphic functions for some complete non-K\"ahler manifolds.

From now on, let $(M, g)$ be a complete Hermitian manifold of complex dimension $n$ and $\omega$ the  associated (1,1)-form. The metirc $g$ is called Gauduchon if $\omega$ satisfies $\partial \bar{\partial} \omega^{n-1}=0$. In \cite{G84}, Gauduchon proved that, when $(M, g)$ is compact, there must exist a  Gauduchon metric $g_{0}$ in the conformal class of $g$. It is interesting to generalize some geometric results from K\"ahler manifolds  to Gauduchon manifolds, for example, the Donaldson-Uhlenbeck-Yau theorem (\cite{NS65}, \cite{DON85}, \cite{UY86}) is  valid for  Gauduchon manifolds (see \cite{Bu, LY, LT,LT95}).

 Let $\varphi $ be a holomorphic function on $M$, it is easy to check that $\partial \overline{\partial } \log |\varphi |\equiv 0$. Following Yau's argument in \cite{Y76}, we consider the following equation
 \begin{equation}\label{log2}
 \widetilde{\Delta}\log u =\sqrt{-1}\Lambda_{\omega }\partial \overline{\partial }  \log u =f,
 \end{equation}
where  $\Lambda_{\omega }$ denotes the contraction with  $\omega $, $u\geq 0$ and the Hausdorff measure of $\{x\in M | u(x)=0\}$ is zero. It is well known that the difference of two Laplacian is given by a first order differential operator as follows
\begin{equation*} \label{laplacian}
(\widetilde{\Delta}-\frac{1}{2}\Delta)\psi= V\cdot \nabla \psi ,
\end{equation*}
where $V$ is a  vector field on $M$. In the non-K\"ahler case, we should handle the first order term, the key is to control the vector field $V$ in order to use the Stokes' theorem which was prove by Yau (Lemma in \cite{Y76}). In fact, we can prove the following theorem.

\begin{thm} \label{theorem1}
Let $(M, g)$ be a complete Gauduchon manifold of complex dimension $n$ with $|d(\omega^{n-1})|\in L^{\infty}(M)$. Suppose $f$ is bounded from below by a constant and is Lebesgue integrable with $0 <\int_{M} f dv_{g} \leq \infty $ or  $f \equiv 0$.  Then there is no non-constant $L^{p}$ smooth solution of the equation (\ref{log2}) for $0<p < \infty$.
\end{thm}

As an application, we obtain the following Liouville theorem.

\begin{thm} \label{theorem2}
Let $(M, g)$ be a complete Gauduchon manifold of complex dimension $n$ with $|d(\omega^{n-1})|\in L^{\infty}(M)$. Then there is no non-constant $L^{p}$ holomorphic functions for $0<p < \infty$.
\end{thm}

{\bf Remark: }{\it Let $(M_{1}, g_{1})$ be a compact Gauduchon manifold and $(M_{2}, g_{2})$ be a complete K\"ahler manifold, it is easy to check that the product Riemannian  manifold $(M_{1}\times M_{2}, g_{1}\times g_{2})$ is a Gauduchon manifold satisfying the assumption in Theorem \ref{theorem1}.}

\medskip

This paper is organized as follows. In Section 2,  we give a proof of Theorems \ref{theorem1}. In Section 3, we consider a vanishing theorem on Higgs bundles over complete non-K\"ahler manifolds.

\section{Proof of Theorem \ref{theorem1}}

In \cite{Y76}, Yau established  the following generalized Stokes theorem which is an extension of Gaffney's result (\cite{Ga}).

\begin{lem} [\protect {\cite[Lemma ]{Y76}}] \label{lemma}
Let $(N, h)$ be a complete Riemannian manifold and $\eta $ be a smooth integrable $(dim N -1)$-form defined on $M$. Then there exists a sequence of domains $B_{i}$ in $N$ such that $N=\cup_{i}B_{i}$, $B_{i}\subset B_{i+1}$ and $\lim_{i\rightarrow +\infty }\int_{B_{i}}d\eta =0$.
\end{lem}

Let $\epsilon>0 $, as that in \cite{Y76}, we set $v_{\epsilon}=(u+\epsilon)^{\frac{1}{2}} $. Form direct computations, it follows that
\begin{equation}
\sqrt{-1}\Lambda_{\omega }\partial \overline{\partial}\log v_{\epsilon}
= \sqrt{-1}\Lambda_{\omega }\partial \overline{\partial}\log u\cdot \frac{u}{2(u+\epsilon)}+\frac{\epsilon|\partial u|^{2}}{2u(u+\epsilon)}
\end{equation}
and
\begin{equation}
\sqrt{-1}\Lambda_{\omega }\partial \overline{\partial}\log v_{\epsilon}
=\sqrt{-1}\Lambda_{\omega }\frac{\partial \overline{\partial} v_{\epsilon}}{v_{\epsilon}}-\sqrt{-1}\Lambda_{\omega }\frac{\partial v_{\epsilon}\wedge \overline{\partial} v_{\epsilon}}{{v_{\epsilon}}^{2}}.
\end{equation}
Thus we have
\begin{equation}\label{awm1}
v_{\epsilon}\sqrt{-1}\Lambda_{\omega }\partial \overline{\partial} v_{\epsilon}
= |\partial v_{\epsilon}|^{2}+\frac{\sqrt{-1}\Lambda_{\omega }\partial \overline{\partial}\log u}{2}\cdot u+\frac{\epsilon|\partial u|^{2}}{2u}
\end{equation}
Let $0<R_1<R_2 $ be any positive numbers. Fix a point $x_{0}\in M$, and choose a nonnegative cut-off function $\varphi $ satisfying
\begin{equation}
\varphi (x)=\left\{ \begin{split}
1, & \quad  x\in B_{x_0}(R_{1}),\\
0, & \quad  x\in M\setminus B_{x_0}(R_{2}),
\end{split} \right.
\end{equation}
$0\leq \varphi \leq 1$ and $|d\varphi|_{\omega }\leq \frac{C}{R_{2}-R_{1}}$,  where $C$ is a positive constant and $B_{x_{0}}(r)$ is the geodesic ball centered at $x_{0}$ with radius $r$.
According to the Stokes formula and the condition $\partial \overline{\partial } (\omega^{n-1})=0$, one can check that
\begin{equation}\label{awm2}
\begin{split}
 & \int_{B(R_2)} {\varphi}^{2}v_{\epsilon}\sqrt{-1}\Lambda_{\omega }\partial \overline{\partial} v_{\epsilon}\frac{{\omega}^{n}}{n!}\\
= & \int_{M} {\varphi}^{2}v_{\epsilon}\sqrt{-1}\partial \overline{\partial} v_{\epsilon} \wedge\frac{{\omega}^{n-1}}{(n-1)!}\\
= & \int_{M} \partial({\varphi}^{2}v_{\epsilon}\sqrt{-1} \overline{\partial} v_{\epsilon} \wedge\frac{{\omega}^{n-1}}{(n-1)!})-\sqrt{-1}\partial {\varphi}^{2}\wedge v_{\epsilon} \overline{\partial} v_{\epsilon} \wedge\frac{{\omega}^{n-1}}{(n-1)!}\\ &-\sqrt{-1}{\varphi}^{2}\partial v_{\epsilon}\wedge \overline{\partial} v_{\epsilon} \wedge\frac{{\omega}^{n-1}}{(n-1)!}+{\varphi}^{2}v_{\epsilon}\sqrt{-1}\overline{\partial}v_{\epsilon}\wedge\partial\frac{{\omega}^{n-1}}{(n-1)!}\\
= & -2\int_{M} \varphi v_\epsilon\langle \partial\varphi, \overline{\partial}v_{\epsilon}\rangle_g \frac{{\omega}^{n}}{n!}-\int_{M}{\varphi}^{2}|\partial v_\epsilon|^{2}\frac{{\omega}^{n}}{n!}+\frac{1}{2}\int_{M} {\varphi}^{2}\sqrt{-1} \overline{\partial} {v_{\epsilon}}^{2}\wedge \partial\frac{{\omega}^{n-1}}{(n-1)!}\\
= & -2\int_{M} \varphi v_\epsilon\langle \partial\varphi, \overline{\partial}v_{\epsilon}\rangle_g \frac{{\omega}^{n}}{n!}-\int_{M}{\varphi}^{2}|\partial v_\epsilon|^{2}\frac{{\omega}^{n}}{n!}\\ &+\frac{1}{2}\int_{M} \sqrt{-1}\overline{\partial}({\varphi}^{2} {v_\epsilon}^{2}\partial \frac{{\omega}^{n-1}}{(n-1)!})-\sqrt{-1}{v_\epsilon}^{2}\overline{\partial}{\varphi}^{2}\wedge\partial\frac{{\omega}^{n-1}}{(n-1)!}\\
\leq &\int_{M}({v_\epsilon}^{2} |\partial \varphi|^{2}+C^\star {v_\epsilon}^{2}\varphi |\partial \varphi|)\frac{{\omega}^{n}}{n!},\\
\end{split}
\end{equation}
where we have used the condition that $|d(\omega^{n-1})|\in L^{\infty}(M)$ in the last inequality .

Combining (\ref{log2}), (\ref{awm1}) and (\ref{awm2}), we obtain
\begin{equation}
\begin{split}
 & \frac{1}{4}\int_{M} {\varphi}^{2}\frac{|\partial u|^{2}}{u+\epsilon}\frac{{\omega}^{n}}{n!}\\
= & \int_{M}{\varphi}^{2} |\partial v_\epsilon|^{2}\frac{{\omega}^{n}}{n!}\\
\leq &\int_{M}{v_\epsilon}^{2}(|\partial \varphi|^{2}+C^\star \varphi|\partial \varphi|)\frac{{\omega}^{n}}{n!}-\int_{M}\frac{f{\varphi}^{2} u}{2}\\
\leq &\ (\frac{C^{2}}{(R_1-R_2)^{2}}+\frac{{C^\star}C}{(R_2-R_1)})\int_{B(R_2)}{v_\epsilon}^{2}\frac{{\omega}^{n}}{n!}-\frac{1}{2}\int_{M} f{\varphi}^{2} u\frac{{\omega}^{n}}{n!},
\end{split}
\end{equation}
and
\begin{equation}
\begin{split}
 & \frac{1}{4}\int_{B(R_2)} {\varphi}^{2}\frac{|\partial u|^{2}}{u}\frac{{\omega}^{n}}{n!}\\
\leq &\ (\frac{C^{2}}{(R_1-R_2)^{2}}+\frac{{C^\star}C}{(R_2-R_1)})\int_{B(R_2)} u^{2}\frac{{\omega}^{n}}{n!}-\frac{1}{2}\int_{B(R_2)} f{\varphi}^{2} u\frac{{\omega}^{n}}{n!}.
\end{split}
\end{equation}
Let $R_1=2R_2 \rightarrow \infty $, we deduce
\begin{equation}
\int_{M} \frac{|d u|^{2}}{u} \frac{{\omega}^{n}}{n!}< \infty,
\end{equation}
and then
\begin{equation}
 \int_{M} |du| \frac{{\omega}^{n}}{n!} \leq \Big(\int_{M} \frac{|d u|^{2}}{u}\frac{{\omega}^{n}}{n!}\Big)^{\frac{1}{2}} \Big(\int_{M} u \frac{{\omega}^{n}}{n!} \Big)^{\frac{1}{2}}< \infty.
\end{equation}

On the other hand, it is easy to see that
\begin{equation}
\sqrt{-1}\Lambda_{\omega }\partial \overline{\partial}\log u_{\epsilon}
= \sqrt{-1}\Lambda_{\omega }\partial \overline{\partial}\log u \cdot\frac{u}{u+\epsilon}+\frac{\epsilon}{u(u+\epsilon)^{2}}|\partial u|^{2}
\end{equation}
and
\begin{equation}\label{Z1}
 \int_{M} |\partial \log u_\epsilon|\frac{{\omega}^{n}}{n!}
=  \int_{M} \frac{|\partial u|}{u+\epsilon}\frac{{\omega}^{n}}{n!}
\leq \int_{M} \frac{|du|}{\epsilon}\frac{{\omega}^{n}}{n!}
 < \infty
\end{equation}
for any $\epsilon >0$, where $u_\epsilon=u+\epsilon $. Using the Stokes formula and the condition $\partial \overline{\partial } (\omega^{n-1})=0$ again, we derive
\begin{equation}
\begin{split}
 &  \sqrt{-1}\Lambda_{\omega }\partial \overline{\partial}\log u_{\epsilon}\frac{{\omega}^{n}}{n!}\\
= & \sqrt{-1}\partial(\overline{\partial} \log u_\epsilon \wedge\frac{{\omega}^{n-1}}{(n-1)!})+\sqrt{-1}\overline{\partial}(\log u_\epsilon\frac{\partial {\omega}^{n-1}}{(n-1)!})\\
= & \sqrt{-1}\partial(\overline{\partial} \log u_\epsilon \wedge\frac{{\omega}^{n-1}}{(n-1)!})+\sqrt{-1}\overline{\partial}(\log (\frac{u}{\epsilon }+1 )\frac{\partial \omega^{n-1}}{(n-1)!}).\\
\end{split}
\end{equation}
Notice that
\begin{equation}\label{Z2}
\int_{M} \log (\frac{u}{\epsilon }+1 ) \frac{\omega^{n}}{n!} \leq \int_{M} \frac{u}{\epsilon }\frac{\omega^{n}}{n!} < \infty .
\end{equation}
 This together with (\ref{Z1}), the condition  that $|d(\omega^{n-1})|\in L^{\infty}(M)$ and lemma \ref{lemma}, gives us that
\begin{equation}\label{awm3}
0=\lim_{i\rightarrow +\infty }\int_{B_{i}} \Delta \log u_{\epsilon} =\int_{M} \frac{fu}{u+\epsilon} \frac{{\omega}^{n}}{n!}+\int_{M} \frac{1}{u(u+\epsilon)^{2}}|d{u}^{2}|\frac{\omega^{n}}{n!}
\end{equation}
for all $\epsilon >0 $, where $B_{i}'s$ are the ones in the lemma \ref{lemma}.
By the same argument  in \cite{Y76}, we know that  $u$ must be a constant. The proof of Theorem \ref{theorem1} is therefore completed.

\section{A vanishing theorem on Higgs bundle}

 Let $(M, \omega )$ be an $n$-dimensional  Hermitian manifold. A Higgs bundle $(E , \overline{\partial }_{E}, \theta )$ over $M$ is a holomorphic bundle $(E , \overline{\partial }_{E})$ coupled with a Higgs field $\theta \in \Omega_M^{1,0}(\mathrm{End}(E))$
such that $\overline{\partial}_{E}\theta =0$ and $\theta \wedge \theta =0$. Higgs bundles first emerged thirty years ago in Hitchin's (\cite{HIT}) study of the self duality
equations on a Riemann surface and in Simpson's subsequent work (\cite{SIM}, \cite{SIM2})  on nonabelian Hodge theory. Such objects have  rich structures and play an important role in many  areas
including gauge theory, K\"ahler and hyperk\"ahler geometry, group representations
and nonabelian Hodge theory. Let $H$ be a Hermitian metric on the bundle $E$, we consider the Hitchin-Simpson connection
 \begin{equation}
 \overline{\partial}_{\theta}:=\overline{\partial}_{E}+\theta , \quad D_{H,  \theta }^{1, 0}:=D_{H }^{1, 0}  +\theta^{*_H}, \quad D_{H,  \theta }= \overline{\partial}_{\theta}+ D_{H,  \theta }^{1, 0},
 \end{equation}
where $D_{H}$ is the Chern connection of $(E,\overline{\partial }_{E},H)$ and $\theta^{*_H}$ is the adjoint of $\theta $ with respect to the metric $H$.
The curvature of this connection is
\begin{equation}
F_{H,\theta}=F_H+[\theta,\theta^{*_H}]+\partial_H\theta+\bar{\partial}_E\theta^{*_H},
\end{equation}
where $F_H$ is the curvature of $D_{H}$ and $\partial_H$ is the $(1, 0)$-part of $D_{H}$.

 Let $s$ be a $\theta $-invariant holomorphic section of a Higgs bundle $(E , \overline{\partial }_{E}, \theta )$, i.e. there exists a holomorphic $1$-form $\eta $ on $M $ such that $\theta (s)=\eta \otimes s$. When the base manifold $(M, \omega )$ is compact, following Kobayashi's techniques \cite{Kobayashi1},  one can obtain vanishing theorems for $\theta $-invariant holomorphic sections on Higgs bundles or Higgs sheaves (see \cite{Bruzzo, Cardo, LZZ}). Now we consider the case that the base manifold is complete non-K\"ahler.

 Since $s$ is $\theta $-invariant, by the formula (4.50) in \cite{LZZ}, we know
  \begin{equation}\label{W01}
\sqrt{-1}\Lambda_{\omega } \langle s, -[\theta , \theta^{\ast H}]s\rangle_{H}=\ |\theta^{\ast H}s- \langle\theta^{\ast H}s,  s\rangle_{H}\frac{s}{|s|_{H}^{2}}|_{H, \omega}^{2}\geq 0.
\end{equation}
A straightforward calculation shows the following Weitzenb\"ock formula
\begin{equation}
\begin{split}
 & \sqrt{-1}\Lambda_{\omega }\partial \overline{\partial}|s|_{H}^{2}\\
= &\ |D_{H}^{1, 0}s|_{H, \omega}^{2}+\sqrt{-1}\Lambda_{\omega }\langle s, F_{H}s\rangle_{H} \\
= &\ |D_{H}^{1, 0}s|_{H, \omega}^{2}-\langle s, \sqrt{-1}\Lambda_{\omega }F_{H, \theta }s\rangle_{H}-\sqrt{-1}\Lambda_{\omega } \langle s, [\theta , \theta^{\ast H}]s\rangle_{H} \\
\geq &\ |D_{H}^{1, 0}s|_{H, \omega}^{2}-\langle s, \sqrt{-1}\Lambda_{\omega }F_{H, \theta }s\rangle_{H}
\end{split}
\end{equation}
on $M$. On the other hand, it holds that
\begin{equation}
\sqrt{-1}\Lambda_{\omega } \partial |s|_{H}^{2}\wedge \overline{\partial }|s|_{H}^{2}\leq |s|_{H}^{2}|D_{H}^{1, 0}s|_{H}^{2}.
\end{equation}
Then these yield that
\begin{equation}\label{log3}
\begin{split}
 & \sqrt{-1}\Lambda_{\omega }\partial \overline{\partial} \log |s|_{H}^{2}\\
= &\frac{1}{|s|_{H}^{2}}\sqrt{-1}\Lambda_{\omega }\partial \overline{\partial}  |s|_{H}^{2}-\frac{1}{|s|_{H}^{4}}\sqrt{-1}\Lambda_{\omega } \partial |s|_{H}^{2}\wedge \overline{\partial }|s|_{H}^{2}\\
\geq &-\frac{1}{|s|_{H}^{2}}\langle s, \sqrt{-1}\Lambda_{\omega }F_{H, \theta }s\rangle_{H}.
\end{split}
\end{equation}
Together with Theorem \ref{theorem1}, we conclude

\begin{thm} \label{theorem3}
Let $(M, \omega )$ be a complete Gauduchon manifold of complex dimension $n$ with $|d(\omega^{n-1})|\in L^{\infty}(M)$, and $(E , \overline{\partial }_{E}, \theta )$ be a Higgs bundle over $M$. If there exists a Hermitian metric $H$ on $E$ such that $\sqrt{-1}\Lambda_{\omega }F_{H, \theta } \leq -f Id_{E}$,  where $f$ is a continuous function on $M$ which is bounded from below and  $\int_{M} f \omega^{n}>0$. Then $E$ admits no non-zero $\theta$-invariant  holomorphic $L^{p}$-sections for $0<p < \infty$.
\end{thm}

\vskip 1 true cm


\bigskip
\bigskip

\noindent {\footnotesize {\it Yuang Li, Chuanjing Zhang and Xi Zhang} \\
{School of Mathematical Sciences, University of Science and Technology of China}\\
{Anhui 230026, P.R. China}\\
{Email:lya1997@mail.ustc.edu.cn; chjzhang@mail.ustc.edu.cn;  mathzx@ustc.edu.cn}

\vskip 0.5 true cm


\begin{thebibliography}{20}

\bibitem{Bruzzo}
U.Bruzzo and B.Gra\~na Otero,
{\em Metrics on semistable and numerically effective Higgs bundles},
J. Reine Angew. Math., {\bf 612}(2007), 59-79.

\bibitem{Bu}
N.P. Buchdahl,
{\em Hermitian-Einstein connections and stable vector bundles over compact complex surfaces},
Math. Ann., {\bf 280}(1988), no. 4, 625-648.


\bibitem{Cardo}
S.A. Cardona,
{\em On vanishing theorems for Higgs bundles},
Differential Geom. Appl., {\bf 35} (2014), 95-102.


\bibitem{DON85}
S.K. Donaldson,
{\em Anti self-dual Yang-Mills connections over complex algebraic surfaces and stable vector bundles},
Proc. London Math. Soc., {\bf50} (1985), no. 1, 1-26.


\bibitem{Ga}
M.P. Gaffney,
{\em A special Stokes's theorem for complete Riemannian manifolds},
Ann. of Math.(2), {\bf 60} (1954), 140-145.


\bibitem{G84}
P. Gauduchon,
{\em Torsion $1$-forms of compact Hermitian manifolds},
Math. Ann., {\bf 267} (1984), no. 4, 495-518.



\bibitem{GW71}
R. E. Greene and H. Wu,
{\em Curvature and complex analysis}, Bull. Amer. Math. Soc., {\bf 77}(1971), 1045-1049.

\bibitem{HIT}
N.J. Hitchin,
{\em The self-duality equations on a Riemann surface},
Proc. London Math. Soc. (3), {\bf55} (1987), no. 1, 59-126.


\bibitem{Kobayashi1}
S. Kobayashi,
{\em Curvature and stability of vector bundles},
Proc. Japan Acad. Ser. A Math. Sci., {\bf 58} (1982), no. 4, 158-162.


\bibitem{LY}
J. Li and S.T. Yau,
{\em Hermitian-Yang-Mills connection on non-K\"ahler manifolds}, Mathematical aspects of string theory (San Diego, Calif., 1986), 560-573, Adv. Ser. Math. Phys., 1, World Sci. Publishing, Singapore, 1987.


 \bibitem{LT95}
M. L\"{u}bke and A. Teleman,
{\em The Kobayashi-Hitchin correspondence}, World
Scientific Publishing Co., Inc., River Edge, NJ, 1995.

\bibitem{LT}
M. L\"{u}bke and A. Teleman,
{\em The universal Kobayashi-Hitchin correspondence on Hermitian manifolds}, Mem. Amer. Math. Soc., {\bf 183} (2006), no. 863.

\bibitem{LZZ}
J.Y. Li, C. J. Zhang, X. Zhang,
{\em Semi-stable Higgs sheaves and Bogomolov type inequality}, Calc. Var. Partial Differential Equations, {\bf 56} (2017), no. 3, Art. 81, 33 pp.


\bibitem{NS65}
M.S. Narasimhan and C.S. Seshadri,
{\em Stable and unitary vector bundles on a compact Riemann surface},
Ann. of Math. (2), {\bf82}(1965), 540-567.


\bibitem{SIM}
C.T. Simpson,
{\em Constructing variations of Hodge structure using Yang-Mills theory and applications to uniformization},
J. Amer. Math. Soc., {\bf1} (1988), no. 4, 867-918.


\bibitem{SIM2}
C.T. Simpson, {\em Higgs bundles and local systems}, Inst. Hautes \'{E}tudes Sci. Publ. Math., {\bf75}(1992), 5-95.


\bibitem{UY86}
K.K. Uhlenbeck and S. T. Yau,
{\em On the existence of Hermitian-Yang-Mills connections in stable vector bundles}, Comm. Pure Appl. Math. {\bf39}(1986), no. S, suppl.,
S257-S293.


\bibitem{Y76}
S.T. Yau,
{\em Some function-theoretic properties of complete Riemannian manifold and their applications to geometry}, Indiana Univ. Math. J., {\bf 25}(1976), no. 7, 659-670.





\end{thebibliography}
\end{document}